%% file: euler4.tex
\documentclass[11pt]{article}
 \input xy
 \xyoption{all}

 \usepackage{esint}

 \usepackage{fullpage}

\newcommand{\rvsR}{rank-varying sub-Riemannian structure}
\newcommand{\ffoot}[1]{}

\newcommand{\segno}{\alpha}
\newcommand{\rotations}{revolutions}
\newcommand{\distribution}{\Upsilon}
\newcommand{\HH}{{\bf (H0)}}
\newcommand{\bD}{\Delta}
\newcommand{\pphi}{\xi}

\renewcommand{\th}{\theta}

\newcommand{\cil}{{\bf W}}
\newcommand{\CIL}{{\bf U}}

\newcommand{\treint}{\sqint}
\newcommand{\f}{f}

\usepackage{mathrsfs}
\usepackage{graphicx}
\usepackage{amssymb}
\usepackage{color}
\newtheorem{theorem}{Theorem}
\newtheorem{corollary}{Corollary}
\newtheorem{lemma}{Lemma}
\newtheorem{proposition}{Proposition}
\newtheorem{definition}{Definition}
\newtheorem{remark}{Remark}

\newcommand{\bt}{\begin{theorem}}
\newcommand{\et}{\end{theorem}}
\newcommand{\bl}{\begin{lemma}}
\newcommand{\el}{\end{lemma}}
\newcommand{\bp}{\begin{proposition}}
\newcommand{\ep}{\end{proposition}}
\newcommand{\bc}{\begin{corollary}}
\newcommand{\ec}{\end{corollary}}
\newcommand{\bdeff}{\begin{definition}}
\newcommand{\edeff}{\end{definition}}
\newcommand{\brem}{\begin{remark}}
\newcommand{\erem}{\end{remark}}
\renewcommand{\r}[1]{(\ref{#1})}
\newcommand{\con}{{\mathcal C}}

\newcommand{\bi}{\begin{itemize}}

\newcommand{\ei}{\end{itemize}}
\newcommand{\bd}{\begin{description}}
\newcommand{\ed}{\end{description}}

\newcommand{\bqn}{\begin{eqnarray}}
\newcommand{\eqn}{\end{eqnarray}}
\newcommand{\eqnn}{\nonumber\end{eqnarray}}

\newcommand{\nn}{\nonumber}
\newcommand{\ba}[1]{\begin{array}{#1}}
\newcommand{\ea}{\end{array}}

\newcommand{\R}{\mathbb{R}}
\newcommand{\Z}{\mathbb{Z}}

\newcommand{\ppotR}[3]
{

\begin{figure}\begin{center}
~\includegraphics[width=#3truecm]{#1.eps}\\
\caption{#2}
\label{#1}
\end{center}
\end{figure}
\noindent$\!\!$}

\newcommand{\g}{\gamma}

\newcommand{\eps}{\varepsilon}

\newcommand{\rvd}{rank-varying distribution}

\newcommand{\VecM}{\mathrm{Vec}(M)}

\newcommand{\Gq}{{\gg}_q}

\newcommand{\Zz}{\mathcal{Z}}
\renewcommand{\gg}{{\bf G}}
\newcommand{\sign}{\mathrm{sign}}
\newcommand{\E}{e}

\title{\LARGE \bf
Two-Dimensional Almost-Riemannian Structures with Tangency Points
}

\author{A.A.~Agrachev\thanks{
 SISSA, via Beirut 2-4, 34014 Trieste, Italy, {\tt agrachev@sissa.it, ghezzi@sissa.it}}
, U.~Boscain\thanks{
Centre de Math\'ematiques Appliqu\'ees, \'Ecole Polytechnique
Route de Saclay, 91128 Palaiseau Cedex, France,
        {\tt ugo.boscain@cmap.polytechnique.fr}}
, G.~Charlot\thanks{
Institut Fourier, UMR 5582, CNRS/Universit\'e Grenoble 1, 100 rue des Maths, BP 74, 38402 St Martin d'H\`eres, France,
{\tt Gregoire.Charlot@ujf-grenoble.fr}} , 
\setcounter{footnote}{1}R.~Ghezzi$^\thefootnote$, \setcounter{footnote}{3}M.~Sigalotti\thanks{INRIA Nancy -- Grand Est, \'Equipe-projet CORIDA, and
Institut \'Elie Cartan, UMR CNRS/INRIA/Nancy Universit\'e, BP 239, 54506 Vand\oe uvre-l\`es-Nancy, France,
        {\tt Mario.Sigalotti@inria.fr}} 
}

\begin{document}
\maketitle
\begin{abstract}
Two-dimensional almost-Riemannian structures are generalized Riemannian structures on surfaces for which a local orthonormal frame is given by a  Lie bracket generating pair of vector fields that can become collinear.
We study the relation between the topological invariants of an almost-Riemannian structure on a compact oriented surface and the rank-two vector bundle over the surface which defines the structure. We analyse the generic case including the presence of tangency points, i.e. points where two generators of the distribution and their Lie bracket are linearly dependent. The main result of the paper provides a classification  of  oriented almost-Riemannian structures on  compact oriented surfaces in terms of  the Euler number of the vector bundle corresponding to the structure. Moreover, we present a Gauss--Bonnet formula for almost-Riemannian structures with tangency points.
\end{abstract}

\section{Introduction}

Let $M$ be a two-dimensional smooth manifold.
A Riemannian distance on $M$ can be seen as the minimum-time function of 
an optimal control problem where admissible velocities are vectors of norm one. 
The control problem can be written locally as 
\bqn
\label{ff}
\dot q=u X(q)+v Y(q)\,,~~~u^2+v^2\leq 1\,,
\eqn
by fixing an orthonormal frame $(X,Y)$. 

Almost-Riemannian structures generalize Riemannian ones by allowing 
$X$ and $Y$ to be collinear at some points. 
If the pair $(X,Y)$ is Lie bracket generating, i.e., if
$$\mathrm{span}\{X(q),Y(q),[X,Y](q),[X,[X,Y]](q),\ldots\}=T_qM$$at every $q\in M$, then~\r{ff}  
is completely controllable and  the
minimum-time function still defines a
continuous distance $d$ on $M$.
Notice that a Riemannian distance can be globally defined on $M$ by a control problem as in~\r{ff} only if the Riemannian structure admits a global orthonormal frame, 
implying that $M$ is parallelizable.
More in general, $X$ and $Y$ are parallel on a  set $\Zz\subset M$ (called {\it singular locus}), which is generically  a one-dimensional embedded submanifold of $M$ (possibly disconnected).

Metric structures
that can be defined {\it locally} by a pair of vector fields
$(X,Y)$ through~\r{ff}
are called {\it almost-Riemannian structures}. An almost-Riemannian structure ${\cal S}$ can be seen as an Euclidean  bundle $E$ of rank two over $M$ (i.e. a vector bundle whose fibre is equipped with a smoothly-varying scalar product $\langle\cdot,\cdot\rangle_q$) and a morphism of vector bundles $\f:E\rightarrow TM$ such that $\f(E_{q})\subseteq T_{q}M$  and the evaluation at $q$ the Lie algebra generated by 
$$
\{\f\circ\sigma\mid\sigma\mathrm{ \, section\,  of\, } E\}
$$
 is equal to $T_{q}M$ for every $q\in M$. 

If $E$ is orientable, we say that ${\mathcal S}$ is {\it orientable}.
If $E$ is isomorphic to the trivial bundle $M\times\R^{2}$, we say that the almost-Riemannian structure is {\it trivializable}. The singular locus $\Zz$  is the set of points $q$  of $M$ at which $\f(E_{q})$ is one-dimensional.
An almost-Riemannian structure is Riemannian if and only if $\Zz=\emptyset$, 
i.e. $\f$ is an isomorphism of vector bundles.

 The first example of
genuinely almost-Riemannian structure
is provided by
the Grushin plane,
which is the almost-Riemannian structure on $M=\R^2$
with $E=\R^2\times\R^2$, $\f((x,y),(a,b))=((x,y),(a,b x))$ and $\langle\cdot,\cdot\rangle$ the canonical Euclidean structure on $\R^2$.
The model was originally introduced in the context of
hypoelliptic operator theory
\cite{FL1,grusin1} (see also \cite{bellaiche,algeria}).
Notice that the
singular locus
is indeed nonempty, being equal to
the $y$-axis.
Another example of (trivializable) almost-Riemannian structure
 appeared in problems of control of quantum mechanical systems
(see \cite{q4,q1}). 

 The notion of almost-Riemannian structure was introduced in \cite{ABS}. In that paper,  an almost-Riemannian structure is  defined  as a locally finitely generated Lie bracket generating
$\con^\infty(M)$-submodule $\bD$ of $\VecM$, the space of smooth vector fields on $M$,
endowed with a bilinear, symmetric map $G:\bD\times\bD\to \con^\infty(M)$ which is positive definite (in a suitable sense). This definition is equivalent to the one given above in terms of Euclidean bundles of rank two over $M$.
A pair of vector fields $(X,Y)$ in $\bD$ is said to be an {\it orthonormal frame for $G$} on some open set $\Omega$ if $G(X,Y)(q)=0$ and
$G(X,X)(q)=G(Y,Y)(q)=1$ for every $q\in\Omega$. Equivalently, there exists a local orthonormal frame $(\sigma,\rho)$ for $\langle\cdot,\cdot\rangle$ such that $X=\f\circ\sigma$, $Y=\f\circ\rho$. 

Almost-Riemannian structures present very interesting phenomena. For instance,  even in the case where
the Gaussian curvature is everywhere negative (where it is defined, i.e., on $M\setminus\Zz$), geodesics may have conjugate points. 
This happens for instance on the Grushin plane (see \cite{ABS} and also \cite{tannaka} in the case of surfaces of revolution.) Moreover it is possible to define non-orientable
almost-Riemannian structures on orientable manifolds and orientable
almost-Riemannian structures on non-orientable manifolds (see \cite{ABS}).

\bigskip
This paper is a continuation of \cite{ABS}, where we provided 
a characterization of generic almost-Riemannian structures
by means of local normal forms, and we  proved a generalization of the Gauss--Bonnet formula. (For generalizations of Gauss--Bonnet formula in related contexts, see \cite{agra-gauss,pelletier}.) Let us briefly recall these results.

The {\it flag} of a submodule $\bD$ of $\VecM$ is the sequence of submodules
$\bD=\bD_1\subset \bD_2\subset\cdots \subset\bD_m \subset \cdots$ 
defined through the recursive formula
$$
\bD_{k+1}=\bD_k+[\bD,\bD_k].
$$
Denote by $\bD_{m}(q)$ the set $\{V(q)\mid V\in \bD_m\}$.
Under generic assumptions, the singular locus $\Zz$ has the following properties: {\bf (i)} $\Zz$ is an
embedded one-dimensional 
submanifold of
$M$;
{\bf (ii)} the points $q\in M$ at which $\bD_2(q)$ is
one-dimensional are isolated;
{\bf (iii)}  $\bD_3(q)=T_qM$ for every $q\in M$.
We  say that $\cal S$ satisfies \HH\ if properties {\bf (i)},{\bf (ii)},{\bf (iii)} hold true. If
 this is the case, a point $q$ of $M$ is called {\it ordinary} if $\bD(q)=T_qM$, {\it Grushin point} if $\bD(q)$ is one-dimensional and $\bD_2(q)=T_qM$, i.e. the distribution is transversal to $\Zz$, and {\it tangency point} if $\bD_2(q)$ is one-dimensional, i.e. the distribution is tangent to $\Zz$. 
If $(\Omega,X,Y)$ is a local generator of $\Delta$  
such that $\Omega\setminus \Zz$ has exactly two components, then 
$(X,Y)$ has different orientations on each of them.   
Normal forms for local generators of $\Delta$  around
ordinary, Grushin and tangency points are recalled in Section~\ref{s-generic}.

The main result of  \cite{ABS} is an extension of the Gauss--Bonnet theorem for orientable almost-Riemannian structures on orientable manifolds, under the hypothesis that there are not tangency points. More precisely, denote by $K:M\setminus\Zz\to\R$ the Gaussian curvature and by $\omega$ a volume form for the Euclidean structure on $E$. Let $dA_s$ be the two-form on $M\setminus \Zz$ given by the pushforward of $\omega$ along $\f$. Fix an orientation $\Xi$ on $M$ and let $M^+$ (respectively, $M^-$) be the subset of $M\setminus\Zz$ where $dA_s$ is a positive (respectively, negative) multiple of $\Xi$.   
The main goal of \cite{ABS} was  to prove the existence and to assign a value to the limit
\bqn
\label{limit}
\lim_{\eps\searrow0}\int_{M_\eps}K\, dA_s,
\eqn
where   $M_\eps=\{q\in M \mid d(q,\Zz)>\eps\}$ and
$d(\cdot,\cdot)$ is the distance globally defined by the almost-Riemannian structure on $M$.
If $M$ has no tangency point, 
then the limit in~\r{limit} can be shown to exist and to be equal to
 $2\pi (\chi(M^+)-\chi(M^-))$, where $\chi$ denotes the Euler characteristic.
When the almost-Riemannian structure
 is trivializable,
 we proved that $\chi(M^+)=\chi(M^-)$ whence the
 limit in~\r{limit}
is equal to zero. Once applied to the special subclass of Riemannian
structures, such a result simply states that the  integral of the
curvature of a parallelizable compact oriented surface (i.e., the
torus) is equal to zero. In a sense, in the standard Riemannian
construction the topology of the surface gives a constraint on the
total curvature through the Gauss--Bonnet formula, whereas for an
almost-Riemannian structure induced by a single pair of vector
fields the total curvature is equal to zero and the topology of the
manifold constrains the metric to be singular on a suitable set.

\bigskip
The main objective of this paper is to complete the analysis in the more complicated case in which ${\cal S}$ has  tangency points.

 The following result provides a classification of almost-Riemannian structures in terms of the Euler number
of the vector bundle associated with it (for a definition of the Euler number see Section~\ref{basdef}).
\bt\label{converse}
Let $M$ be a compact
oriented
two-dimensional manifold  endowed with an oriented almost-Riemannian structure
 ${\cal S}=(E,\f,\langle\cdot,\cdot\rangle)$ satisfying the generic hypothesis \HH.
Then  $\chi(M^+)-\chi(M^-)+\tau( {\cal S})=\E(E)$, where  
$\E(E)$ denotes the Euler number of  $E$ and 
$\tau({ \mathcal S})$ is the number of \rotations\ of 
$\Delta$ on $\Zz$ computed with respect to  the orientation induced by $M$ on $\Zz$. 
\et
For a detailed definition of $\tau({\cal S})$ see Section~\ref{deftau}.
Notice that the Euler number $\E(E)$  measures  how far the vector bundle $E$ is from 
the trivial one. Indeed,  $E$ is trivial if and only if $\E(E)=0$. As a direct consequence of Theorem~\ref{converse} 
we get that ${\cal S}$ is trivializable if and only if $\chi(M^+)-\chi(M^-)+\tau({\cal S})=0$.

\bigskip
For what concerns the notion of integrability of the curvature with respect to the Riemannian density on $M\setminus \Zz$, it turns out that the hypothesis made in \cite{ABS} about the absence of tangency points is not just technical. Indeed, in Section~\ref{numsim} we present  some numerical simulations  strongly hinting that the limit in~\r{limit} diverges, in general, if tangency points are present. One possible explanation of this fact is the interaction between different orders in the asymptotic expansion of the almost-Riemannian distance.
In this paper, to avoid this interference,  we define a 3-scale integral of the curvature. Let ${\cal T}$ be the set of tangency points of ${\cal S}$. Associate with every $q\in {\cal T}$, 
a two parameters ``rectangular''   neighborhood $B^q_{\delta_1,\delta_2}$ ($\delta_1$ and $\delta_2$ 
playing the role of  lengths of the sides of the rectangle) built as follows.
We consider a smooth curve 
$(-1,1)\ni s\mapsto w(s)$
passing through the tangency point $w(0)=q$ and transversal to the distribution at $q$. We then consider, for each $s\in (-1,1)$, 
 the  geodesic $\g_s$ (parameterized by arclength) 
 such that $\g_s(0)=w(s)$ and minimizing locally the distance from $\{w(s)\mid s \in(-1,1) \}$ (such geodesic exists and is unique due to the transversality of the curve $w$; for details see \cite{ABS}). 
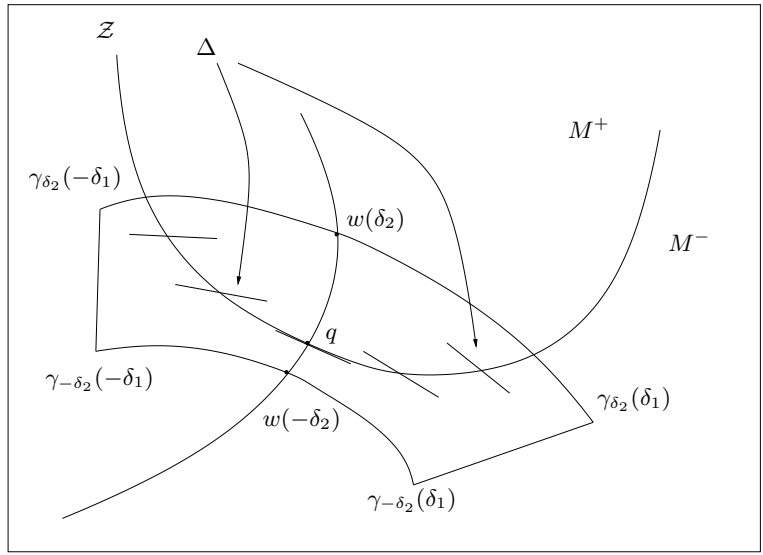
\begin{figure}[h!]
\begin{center}
\input{boxtrasv.pstex_t}
\caption{The rectangular box $B^q_{\delta_1,\delta_2}$}\label{boxbr}
\end{center}
\end{figure}
For $\delta_1,\delta_2$ sufficiently small, the rectangle $B^q_{\delta_1,\delta_2}$ is the subset of $M$ 
 containing the tangency point $q$
 and having as boundary 
$$\g_{\delta_2}([-\delta_1,\delta_1] )\cup
\g_{[-\delta_2,\delta_2]}(\delta_1)\cup
\g_{-\delta_2}([-\delta_1,\delta_1])\cup
\g_{[-\delta_2,\delta_2]}(-\delta_1)
$$
(see Figure~\ref{boxbr}). Let $M_{\eps,\delta_1,\delta_2}=M_\eps \setminus  \bigcup_{q\in {\cal T} }  B^q_{\delta_1,\delta_2}$.
We say that $K$ is {\it 3-scale ${\cal S}$-integrable} if 
\bqn
\lim_{\delta_1\to0}\lim_{\delta_2\to0}\lim_{\eps\to0}\int_{M_{\eps,\delta_1,\delta_2}} K dA_s
\eqn
exists and is finite. In this  case we denote such a limit by $\treint_M K dA_s$. 
The following result, proved in Section~\ref{secpro}, is a generalization of the classical Gauss--Bonnet
 formula for Riemannian structures
to generic oriented two-dimensional almost-Riemannian structures. It concludes the analysis started in \cite{ABS}
including the presence of tangency points.
\bt\label{t-GB}
Let $M$ be a compact
oriented
two-dimensional manifold. For an oriented almost-Riemannian structure
on $M$ satisfying the generic hypothesis \HH, $K$ is 3-scale ${\mathcal S}$-integrable and
\bqn
\treint_M K dA_s=2\pi\E(E).
\eqnn
\et
Notice that the construction of $M_{\eps,\delta_1,\delta_2}$ is not intrinsic since it depends on the choice of a manifold transversal to $\Zz$ at each tangency point and on its parameterization. The result, however, is.
An interesting question is whether a canonical way of choosing these manifolds and their parameterizations exists. 
This is related to the problem of finding intrinsic local normal forms for almost-Riemannian structures.

As a consequence  of Theorems~\ref{converse} and~\ref{t-GB} we get the following corollary.
\bc
Let $M$ be a compact
oriented
two-dimensional manifold without boundary. For an oriented almost-Riemannian structure ${\cal S}$
on $M$ satisfying the generic hypothesis \HH\ we have 
$\treint_M K dA_s=0$ if and only if ${\cal S}$ is trivializable. 
In particular, if $\cal S$ has not tangency points then
$\int_M K dA_s=0$ if and only if ${\cal S}$ is trivializable. 

\ec
These results complete the analysis of the relation between the integral of the curvature and the topology of the manifold for two-dimensional almost-Riemannian structures.

\bigskip
The structure of the paper is the following. In Section~\ref{basdef} we introduce the notion of $n$-dimensional rank-varying sub-Riemannian structure and recall some definitions and results given in \cite{ABS}. In Section~\ref{deftau}
 we define the number of \rotations\ of a one-dimensional distribution along an oriented closed curve. 
In Section~\ref{proof} we prove Theorem~\ref{converse}. First we construct a section of the sphere bundle on a tubular
neighborhood of the singular locus having each tangency point as singularity. Then, we extend the section to the 
complement in the manifold of a finite set and show the sum of the indeces at all the singularities to be equal to
 $\chi(M^+)-\chi(M^-)+\tau({\cal S})$.
 In Section~\ref{integr}, the relation between tangency points and integrability of the curvature with respect to the area form associated with an almost-Riemannian structure is discussed. In particular, we provide in Section~\ref{numsim} some numerical simulations strongly hinting that  in presence of tangency points the integral $\int_{M_{\eps}}KdA_s$ does not converge as $\eps$ tends to zero. This leads us to introduce in Section~\ref{3int} the notion of 3-scale ${\cal S}$-integrability. Thanks to a Gauss--Bonnet formula for almost-Riemannian surfaces with boundary given in \cite{high-order}, we compute in Section~\ref{secpro} the total curvature of a generic two-dimensional almost-Riemannian structure with tangency points, proving Theorem~\ref{t-GB}.

\section{Basic definitions}\label{basdef}


Let $M$ be a   $n$-dimensional manifold. Throughout the paper, unless specified,  manifolds are smooth (i.e. $\con^{\infty}$) and without boundary; vector fields  and differential forms are smooth. Given a vector bundle $E$ over $M$, the set of smooth sections of $E$, denoted by $\Gamma(E)$,
is naturally endowed with the structure of $\con^\infty(M)$-module. In the case $E=TM$ we denote $\Gamma(E)$ by $\VecM$. 

Given  an oriented vector bundle of rank $n$ over a compact connected oriented $n$-manifold $M$,  the Euler number of $E$, denoted by $\E(E)$, is the self-intersection number of $M$ in $E$, where $M$ is identified with the zero section. To compute $\E(E)$, consider a smooth section $\sigma:M\rightarrow E$ transverse to the zero section. Then, by definition,
$$
\E(E)=\sum_{p\mid \sigma(p)=0}i(p,\sigma),
$$
where $i(p,\sigma)=1$, respectively $-1$, if $d_p\sigma:T_pM\rightarrow T_{\sigma(p)}E$ preserves, respectively reverses, the orientation. Notice that if we reverse the orientation on $M$ or on $E$  then $\E(E)$ changes sign. Hence, the Euler number of an orientable vector bundle $E$ is defined up to a sign, depending on the orientations of both $E$ and $M$. Since reversing the orientation on $M$ also reverses the orientation of $TM$, the Euler number of $TM$ is defined unambiguously and is equal to $\chi(M)$, the Euler characteristic of $M$.
\begin{remark}
 Assume that $\sigma\in\Gamma(E)$ has only isolated zeros, i.e. the set $\{p\mid\sigma(p)=0\}$ is finite. If $E$ is endowed with a smooth scalar product $\langle\cdot,\cdot\rangle$, we define $\tilde{\sigma}:M\setminus\{p\mid\sigma(p)=0\}\rightarrow SE$ by $\tilde{\sigma}(q)=\frac{\sigma(q)}{\sqrt{\langle\sigma(q),\sigma(q)\rangle}}$. Then if $\sigma(p)=0$, $i(p,\tilde{\sigma})=i(p,\sigma)$ is equal to the degree of the map $\partial B\rightarrow S^{n-1}$ that associate with each $q\in\partial B$ the value $\tilde{\sigma}(q)$, where $B$ is a neighborhood of $p$ diffeomorphic to an open ball in $\R^n$ that does not contain any other zero of $\sigma$.

Notice that if $i(p,\sigma)\neq 0$, the limit $\lim_{q\rightarrow p}\tilde{\sigma}(q)$ does not exist and, in this case, we say that $\tilde{\sigma}$ has a singularity at $p$.  
\end{remark}

\begin{definition}\label{fiberrvd}
A {\it $(n,k)$-\rvd} on a $n$-dimensional manifold $M$ is 
a pair $(E,\f)$ where $E$ is a vector bundle of rank $k$ over $M$ and $\f:E\rightarrow TM$ is a morphism of vector bundles satisfying {\bf (i)} $\f$ induces the identity on $M$ i.e. $\f(E_q)\subseteq T_qM$ for every point $q\in M$; {\bf (ii)} the map $\sigma\mapsto\f\circ\sigma$ from $\Gamma(E)$ to $\VecM$ is injective.   
\end{definition}

\brem\label{equiv} 
Definition~\ref{fiberrvd} is equivalent to the definition of $(n,k)$-rank varying distribution given in \cite{ABS}, namely a submodule $\bD$ of $\VecM$ locally generated by $k$ (and not less than $k$) vector fields. Given a $(n,k)$-\rvd\ $(E,\f)$, we define $\bD=\{\f\circ\sigma\mid\sigma\in\Gamma(E)\}$. Thanks to {\bf (i)} $\bD$ is a submodule of $\VecM$. The condition {\bf (ii)} implies that $\bD$ cannot be generated by less than $k$ vector fields. In the following, we use either definition, depending on our convenience.
\erem

Let $(E,\f)$ be a $(n,k)$-\rvd, $\bD=\{\f\circ\sigma\mid\sigma\in\Gamma(E)\}$ be its associated submodule and denote by $\bD(q)$  the linear subspace $\{V(q)\mid  V\in \bD\}=\f(E_q)\subseteq  T_q M$.
Let  $\mathrm{Lie}(\bD)$ be the smallest Lie subalgebra of
$\mathrm{Vec}(M)$
containing $\bD$ and $\mathrm{Lie}_q(\bD)=\{V(q)\mid V\in \mathrm{Lie}(\bD)\}$ for
every $q\in M$.
We say that $(E,\f)$ satisfies the {\it Lie bracket generating condition} if
$\mathrm{Lie}_q(\bD)=T_q M$ for every $q\in M$.

A property $(P)$ defined for $(n,k)$-\rvd s is said to be {\it generic}
if for every vector bundle $E$ of rank $k$ over $M$, $(P)$ holds for every $\f$ in an  open and dense  subset of the set of  morphisms of vector bundles from $E$ to $TM$ inducing the identity on $M$, endowed with the 
$\con^\infty$-Whitney topology. 
E.g., generically, a $(n,k)$-\rvd\ is  Lie bracket generating, provided that $k>1$.

We say that a $(n,k)$-\rvd\ $(E,\f)$ is {\it orientable} if $E$ is orientable as a vector bundle.
 Similarly, $(E,\f)$ is {\it trivializable} if $E$ is isomorphic to the trivial bundle $M\times\R^k$. In this case, by definition, we have $\E(E)=0$. 

A rank-varying sub-Riemannian structure is defined by requiring that $E$ is an Euclidean bundle.

\begin{definition}\label{gensrs}
A {\it $(n,k)$-\rvsR}  is a triple ${\cal S}=(E,\f,\langle\cdot,\cdot\rangle)$ where $(E,\f)$ is a Lie bracket generating $(n,k)$-\rvd\ on a  manifold $M$ and $\langle\cdot,\cdot\rangle_q$ is a   scalar product on $E_q$ smoothly depending on $q$. 
\end{definition}

\brem\label{equival} Definition~\ref{gensrs}  is equivalent to Definition 4 in \cite{ABS}, with 
$G:\bD\times\bD\to\con^\infty(M)$ defined by $G(V,W)=\langle\sigma_V,\sigma_W\rangle$ where $\sigma_V,\sigma_W$ are the unique sections of $E$ satisfying $\f\circ\sigma_V=V,\f\circ\sigma_W=W$.
\erem

\bigskip
Definition~\ref{gensrs} generalizes several classical structures. First of all, a Riemannian structure $(M,g)$ is a $(n,n)$-\rvsR\ with $E=TM$, $\f=1_{TM}$ and $\langle\cdot,\cdot\rangle=g(\cdot,\cdot)$. Classical sub-Riemannian structures (see \cite{book2,montgomery}) are $(n,k)$-\rvsR s such that $E$ is a proper Euclidean subbundle of $TM$ and $\f$ is the inclusion.
Finally, we define {\it $n$-dimensional almost-Riemannian structures} ($n$-ARSs for short) as $(n,n)$-\rvsR s.
\bigskip


Let ${\cal S}=(E,\f,\langle\cdot,\cdot\rangle)$ be a  $(n,k)$-\rvsR. 
For every $q\in M$ 
and every $v\in\bD(q)$ define
\bqn
\Gq(v)=\inf\{\langle u, u\rangle_q \mid u\in E_q,\f(u)=v\}.
\eqnn

If $\sigma_1,\dots,\sigma_k$ is an orthonormal frame for $\langle\cdot,\cdot\rangle$ on an open subset $\Omega$ of $M$, an {\it  orthonormal frame for $G$} on $\Omega$ is given by  $\f\circ\sigma_1,\dots,\f\circ\sigma_k$. One easily proves that 
orthonormal frames are systems of local generators of $\bD$. Notice that ${\cal S}$ is trivializable if and only if there exists a global orthonormal frame for $\langle\cdot,\cdot\rangle$. Hence ${\cal S}$ is trivializable if and only if there exists a global orthonormal frame for $G$.

A curve $\g:[0,T]\to M$ is said to be {\it admissible} for ${\cal S}$ 
if  it is Lipschitz continuous and  there exists a measurable essentially bounded function 
\bqn
[0,T]\ni t\mapsto u(t)\in E_{\g(t)}
\eqnn
called  {\it control function}, such that 
$\dot \g(t)=\f(u(t))$  for almost every $t\in[0,T]$.
Given an admissible 
curve $\g:[0,T]\to M$, the {\it length of $\g$} is  
\bqn
\ell(\g)= \int_{0}^{T} \sqrt{ \gg_{\gamma(t)}(\dot \g(t))}~dt.\eqnn
The function $\ell(\g)$ is invariant under reparameterization of the curve $\g$. Moreover, if an admissible curve $\g$ minimizes the {\it energy functional} 
$
J(\g)=\int_0^T \gg_{\gamma(t)}(\dot \g(t))~dt
$
(with fixed $T$ and fixed endpoints) 
then $v=\sqrt{\gg_{\gamma(t)}(\dot \g(t))}$ is constant and 
$\g$ is also a minimizer of $\ell(\cdot)$. 
On the other hand a minimizer $\g$ of $\ell(\cdot)$ such that  $v$ is constant is a minimizer of $J(\cdot)$ with $T=\ell(\g)/v$.
A {\it geodesic} for  ${\cal S}$  is an admissible 
curve $\g:[0,T]\to M$ such that 
for every sufficiently small interval 
$[t_1,t_2]\subset [0,T]$, $\g|_{[t_1,t_2]}$ is a minimizer of $J(\cdot)$. 
A geodesic for which $\gg_{\gamma(t)}(\dot \g(t))$ is (constantly) 
equal to one is said to be parameterized by arclength. 

The  distance induced by ${\cal S}$ on $M$ is defined as
\bqn\label{e-dipoi}
d(q_0,q_1)=\inf \{\ell(\g)\mid \g(0)=q_0,\g(T)=q_1, \g\ \mathrm{admissible}\}.
\eqn
The finiteness and the continuity of $d(\cdot,\cdot)$ with respect 
to the topology of $M$ are guaranteed by  the Lie bracket generating 
assumption on the \rvsR\ (see \cite{velimir}).  
The distance $d(\cdot,\cdot)$ endows $M$ with the 
structure of metric space compatible with the topology of $M$ as smooth manifold.

\subsection{Normal forms for generic 2-ARSs}\label{s-generic}


Given a 2-ARS ${\mathcal
S}$, its {\it singular locus} is the set of points $q$ where $\bD(q)=\f(E_q)$  has not maximal rank, that is,
$$\Zz=\{q\in M\mid \mbox{dim}(\Delta(q))=1\}.$$ 
Notice that $\mbox{dim}(\Delta(q)) > 0$ for every point $q \in M$, because $\Delta$ is a bracket generating distribution.\\ 

We say that $\cal S$ {\it satisfies condition} \HH\  if the following properties hold:
{\bf (i)} $\Zz$ is an
embedded one-dimensional 
submanifold of
$M$;
{\bf (ii)} the points $q\in M$ at which $\bD_2(q)$ is
one-dimensional are isolated;
{\bf (iii)}  $\bD_3(q)=T_qM$ for every $q\in M$, where $\bD_1 = \bD$ and $\bD_{k+1}=\bD_k+[\bD,\bD_k]$.

\begin{proposition}[\cite{ABS}] \label{p-generic}
Property \HH\ is generic for 2-ARSs.
\end{proposition}

ARSs satisfying hypothesis \HH\ admit the following local normal forms. 
\begin{theorem}[\cite{ABS}]
\label{t-normal}
Given a  2-ARS ${\mathcal S}$ satisfiyng \HH, for every point
$q\in M$ there exist a neighborhood $U$ of $q$ and an orthonormal frame
$(X,Y)$ for $G$ on $U$, such that up
to a smooth change of coordinates defined on $U$, $q=(0,0)$ and $(X,Y)$
has one of the
forms
\bqn
\mathrm{(F1)}&& ~~X(x,y)=(1,0),~~~Y(x,y)=(0,e^{\phi(x,y)}), \nn  \\
\mathrm{(F2)}&& ~~X(x,y)=(1,0),~~~Y(x,y)=(0,x e^{\phi(x,y)}),\nn   \\
\mathrm{(F3)}&& ~~X(x,y)=(1,0),~~~Y(x,y)=(0,(y -x^2
\psi(x))e^{\pphi(x,y)}), \nn
\eqn
where $\phi$, $\pphi$ and $\psi$ are smooth real-valued functions such that
$\phi(0,y)=0$ and  $\psi(0)>0$.
\et
Let ${\mathcal S}$ be a 2-ARS satisfying \HH.
 A point $q\in M$ is said to be an
{\it ordinary point} if $\bD(q)=T_q M$, hence, if ${\mathcal
S}$ is locally described by (F1). We call $q$ a  {\it Grushin
point} if $\bD(q)$ is one-dimensional and $\bD_2(q)=T_q M$, i.e. if
the local description (F2) applies. Finally, if
$\Delta(q)=\Delta_2(q)$ has dimension one and $\bD_3(q)=T_q M$
then we say that $q$ is a {\it tangency point} and ${\mathcal
S}$ can be described near $q$ by the normal form (F3). We define
 $${\cal T}=\{q\in \Zz\mid q \mbox{ tangency point of } {\cal S}\}.$$



\subsection{A Gauss--Bonnet formula for 2-ARSs}

Let $M$ be an orientable two-dimensional manifold and ${\cal S}=(E,\f,\langle\cdot,\cdot\rangle)$ an orientable 2-ARS on $M$.
Notice that $\langle\cdot,\cdot\rangle$ defines a Riemannian structure on $M\setminus \Zz$. Denote by $K$ the Gaussian curvature of such a structure and by $\omega$ a volume form for the Euclidean structure on $E$. Let $dA_s$ be the two-form on $M\setminus \Zz$ given by the pushforward of $\omega$ along $\f$. 
Once an orientation on $M$ is fixed, $M\setminus \Zz$ is split into two open sets
$M^+$ and $M^-$ such that 
the orientation on $M$ coincides with the one defined by $dA_s$ on $M^+$
and with its opposite on $M^-$. 

For every $\eps>0$ let $M_\eps=\{q\in M\mid d(q,\Zz)>\eps\}$, where  $d(\cdot,\cdot)$ is 
the almost-Riemannian distance (see equation~\r{e-dipoi}).  
We say that $K$ is {\it ${\cal S}$-integrable} if 
\bqn
\lim_{\eps\to0}\int_{M_\eps}K~dA_s
\eqnn
exists and is finite. In this  case we denote such limit by $\int_M K dA_s$.

When ${\cal S}$ has no tangency points $K$ happens to be  ${\cal S}$-integrable
and $\int_M K dA_s$ is determined by the topology of $M^+$ and $M^-$. This result, recalled in Theorem~\ref{gbsm}, can be seen as a
generalization of Gauss--Bonnet formula to ARSs. 

\bt[\cite{ABS}]\label{gbsm}
Let $M$ be a compact
oriented
two-dimensional manifold, endowed with an oriented 2-ARS ${\cal S}$ for which condition \HH\ holds true.  
Assume that ${\cal S}$ has no tangency points. 
Then  $K$ is ${\cal S}$-integrable and $\int_M K dA_s=2\pi(\chi(M^+)-\chi(M^-))$.
\et

\section{Number of \rotations\ of $\bD$}\label{deftau}

Let $M$ be a compact
oriented
two-dimensional manifold and $W$ an oriented closed simple curve in $M$. Since $M$ is oriented, $TM|_{W}$ is isomorphic to the trivial bundle $W\times\R^2$ and its projectivization is isomorphic to $W\times S^1$. Every subbundle $\distribution$ of $TM|_W$ of rank one can be seen as a section of the projectivization of $TM|_W$ i.e. a smooth map $\distribution:W\rightarrow W\times S^1$ such that $\pi_1\circ\distribution=\mathrm{Id}_W$, where $\pi_1:W\times S^1\rightarrow W$ denotes the projection on the first component. We define $\tau(\distribution,W)$, the {\it number of \rotations\ } of $\distribution$ along $W$, to be the degree of the map $\pi_2\circ\distribution:W\rightarrow S^1$, where $\pi_2:W\times S^1\rightarrow S^1$ is the projection on the second component. Notice that $\tau(\distribution,W)$ changes sign if we reverse the orientation of $W$.

\begin{remark}\label{fico} Let $M,W,\distribution$ as above and let us show how to compute $\tau(\distribution,W)$. Let $V\in\Gamma(TW)$ be a never-vanishing vector field along $W$. Then $\mathrm{span}(V)$ is a subbundle of $TM|_W$ of rank one and there exists an isomorphism $t:TM|_W\rightarrow W\times\R^2$ such that $t\circ V(q)=(q,(1,0))$. The trivialization $t$ induces an isomorphism, still denoted by $t$, between  the projectivization of $TM|_W$ and $W\times S^1$ such that
$\pi_2\circ t\circ V:W\rightarrow S^1$ is constant. To simplify notations, we omit $t$ in the following. Let $\pi_2\circ V(q)\equiv \th_0$. Assume that $\th_0$ is  a regular value of $\pi_2\circ \distribution$ (otherwise there exists a smooth section $\tilde{\distribution}$ homotopic to $\distribution$ having $\th_0$ as a regular value). 
 By definition,
$$
\tau(\distribution,W)=\sum_{q\mid\pi_2\circ\distribution(q)=\th_0}\sign( d_q(\pi_2\circ\distribution)),
$$
where $d_q$ denotes the differential at $q$ of a smooth map. Notice that a  point $q$ satisfies $\pi_2\circ\distribution(q)=\th_0$ if and only if $\distribution(q)$ is tangent to $W$ at $q$.
\end{remark}

Let $M$ be a compact
oriented
two-dimensional manifold and
${\cal S}$ be an oriented 2-ARS satisfying \HH. Then the singular locus $\Zz$ is the boundary of $M^+$. Fix on $\Zz$ the orientation induced by $M^+$ and let $C(\Zz)$ denote the set of connected component of $\Zz$.
Since $\bD$ is one-dimensional along $\Zz$, we can define
$\tau({\cal S})=\sum_{W\in C(\Zz)}\tau(\bD,W)$.

\brem\label{nuovopn}
Let ${\cal S}$ be a 2-ARS satisfying hypothesis \HH. Let $V,\th_0,\pi_2$ be as in Remark~\ref{fico} with $\distribution=\bD$ and $W\in C(\Zz)$. Since $\bD_3(q)=T_qM$ for every point $q\in M$,  $\th_0$ is a regular value of $\pi_2\circ\bD$. Moreover, the set of points $q\in W$ such that $\pi_2\circ\bD(q)=\th_0$ is the set of tangency points of ${\cal S}$ belonging to $W$. Hence,
$$
\tau(\bD,W)=\sum_{q\in W\cap {\cal T}}\sign(d_q(\pi_2\circ\bD)).
$$
We define $\tau_q=\sign(d_q\pi_2\circ\bD)$ (see Figure \ref{tauu}). Clearly, $\tau({\cal S})=\sum_{q\in{\cal T}}\tau_q$.
\erem

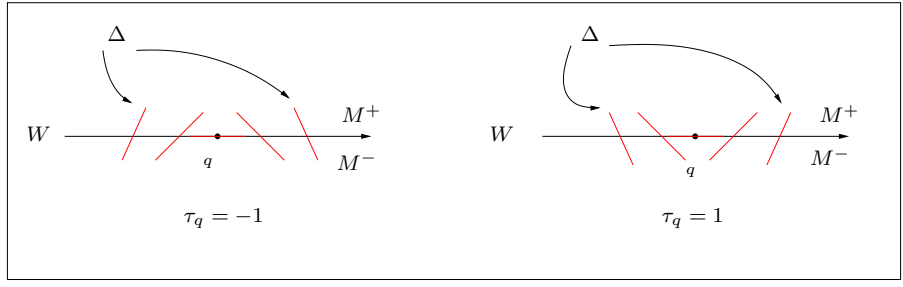
\begin{figure}[h!]
\begin{center}
\input{tauq.pstex_t}

\caption{Tangency points with opposite contributions}\label{tauu}
\end{center}
\end{figure}

\section{Proof of Theorem~\ref{converse}}\label{proof}

The idea of the proof is to find a section $\sigma$ of $SE$ with isolated singularities $p_1,\dots,p_m$ such that $\sum_{j=1}^mi(p_j,\sigma)=\chi(M^+)-\chi(M^-)+\tau({\cal S})$. 
In the sequel, we consider $\Zz$ to be oriented with the orientation induced by $M^+$.

We start by defining $\sigma$ on a neighborhood of $\Zz$. Let $W$ be a connected component of $\Zz$. Since 
$M$ is oriented, there exists an open tubular neighborhood $\cil$ of $W$ and a diffeomorphism
$\Psi:S^1\times (-1,1)\rightarrow \cil$ that preserves the orientation and $\Psi|_{S^1\times\{0\}}$ is an orientation-preserving diffeomorphism between $S^1$ and $W$. 
%
%
 Remark that  $\left.\f:E\right|_{\cil^+}\rightarrow T\cil^+$ is an orientation-preserving isomorphism of vector bundles, while $\left.\f:E\right|_{\cil^-}\rightarrow T\cil^-$ is an orientation-reversing isomorphism of vector bundles, where $\cil^\pm=\cil\cap M^\pm$. For every $s\neq 0$,  lift  the tangent vector to $\th\mapsto\Psi(\th,s)$  to $E$ using $\f^{-1}$,   rotate it  by the angle $\pi/2$ and normalize it: $\sigma$ is defined as this unit vector (belonging to $E_{\Psi(\th,s)}$) if $s>0$, its opposite if $s<0$.  In other words, $\sigma:\cil\setminus W\rightarrow SE$ is given by  
\begin{equation}\label{sez}
\sigma(q)= \sign (s)\frac{R_{\pi/2}\f^{-1}(\frac{\partial \Psi}{\partial \th}(\th,s))}{\sqrt{\langle \f^{-1}(\frac{\partial \Psi}{\partial \th}(\th,s)), \f^{-1}(\frac{\partial \Psi}{\partial \th}(\th,s))\rangle}}, \,\,\,(\th,s)=\Psi^{-1}(q),
\end{equation}
where $R_{\pi/2}$ denotes the rotation (with respect to the Euclidean structure) in $E$ by angle $\pi/2$ in the counterclockwise sense.   
The following lemma shows that $\sigma$ can be extended to a continuous section from
$\cil\setminus {\cal T}$ to $SE$.

\bl\label{unico}
 $\sigma$ can be continuously extended to every point $q\in W\setminus {\cal T}$.
\el

{\bf Proof.}
Let $q \in W\setminus{\cal T}$, $U$ be a neighborhood of $q$ in $M$ and $(x,y)$ be a system of coordinates on $U$ centered at $q$ such the almost-Riemannian structure has the form (F2) (see Theorem~\ref{t-normal}). Assume, moreover, that $U$ is a trivializing neighborhood of both $E$ and $TM$ and the pair of vector fields $(X,Y)$ is the image under $\f$ of a positively-oriented local orthonormal frame of $E$.
Then $W \cap U =\{(x,y)\mid x=0\}$. Since $\frac{\partial \Psi}{\partial \th}(\th,0)$
 is non-zero and tangent to $W$,    $\frac{\partial \Psi}{\partial \th}(\th,0)$ is tangent to the  $y$-axis.
Hence, thanks to the Preparation Theorem  \cite{malgrange}, there exist $h_2:\R\rightarrow\R$, $h_1,h_3:\R^2\rightarrow\R$ smooth functions such that $h_2(y)\neq 0$ for every $y\in \R$ and for $\Psi(\th,s)\in  U$
$$
\frac{\partial\Psi}{\partial \th} (\th,s) = (x h_1(x,y), h_2(y)+x h_3(x,y)),
$$
where $(x,y)$ are the coordinates of the point $\Psi(\th,s)$. Let us compute $\sigma$ at a point $p\in (\cil\cap U)\setminus W$.
Since 
$$
\frac{\partial\Psi}{\partial \th} (\th,s)=x h_1(x,y)X(x,y)+\frac{h_2(y)+xh_3(x,y)}{x e^{\phi(x,y)}}Y(x,y), 
$$
then
$$
\f^{-1}\left(\frac{\partial\Psi}{\partial \th} (\th,s)\right)=x h_1(x,y)\zeta(x,y)+\frac{h_2(y)+xh_3(x,y)}{x e^{\phi(x,y)}}\rho(x,y), 
$$
where $(\zeta,\rho)$ is the unique local orthonormal basis of $E|_U$ such that $\f\circ\zeta=X$ and $\f\circ\rho=Y$.
Notice that  $U\cap M^+=\{(x,y)\mid x>0\}$ and $U\cap M^-=\{(x,y)\mid x<0\}$.
Using formula~\r{sez}, for $(x,y)=\Psi(\th,s)\in U\setminus W$ one easily gets
$$
\sigma(x,y)=\frac{\sign (x)}{l(x,y)}\left( -\frac{h_2(y)+x h_3(x,y)}{x e^{\phi(x,y)}}\zeta+x h_1(x,y)\rho\right),
$$
where 
$$
l(x,y)=\sqrt{ x^2h_1(x,y)^2+\frac{(h_2(y)+xh_3(x,y))^2}{x^2e^{2\phi(x,y)}}}.
$$
Since
$$
\lim_{x\rightarrow 0}\frac{\sign (x) (h_2(y)+x h_3(x,y))}{l(x,y)x e^{\phi(x,y)}}=\frac{h_2(y)}{|h_2(y)|}~~~~\mathrm{ and }~~~\lim_{x\rightarrow 0}\frac{\sign (x) x h_1(x,y)}{l(x,y)}=0,
$$
$\sigma$ can be continuously extended to the set $\{x=0\}=W\cap U$. 
\hfill$\blacksquare$

The next step of the proof is to show that for every $q\in W\cap {\cal T}$, $i(q,\sigma)=\tau_q$. 
 
\begin{lemma}
Let $\sigma:\cil\setminus {\cal T}\rightarrow SE$ be the continuous section obtained in Lemma~\ref{unico}. Then, for every 
$q\in W \cap {\cal T}$
the index of $\sigma$ at $q$ is equal to $\tau_q$ and, consequently,
\begin{equation}\label{verita}
\sum_{q\in {\cal T}\cap W}i(q,\sigma)=\tau(\bD,W).
\end{equation}
\end{lemma}

{\bf Proof.}
Let $q \in W\cap{\cal T}$, $U$ be a neighborhood of $q$ in $M$ and $(x,y)$ be a system of coordinates  on $U$ centered at $q$ such the almost-Riemannian structure has the form (F3) (see Theorem~\ref{t-normal}) i.e.
a local orthonormal frame $(X,Y)$ is given by
$$
X=(1,0), ~~~~~Y=(0,(y-x^2\psi(x))e^{\xi(x,y)}).
$$
Define $\segno=1$, respectively $-1$, if $(X,Y)$ is the image under $\f$ of a positively-oriented, respectively negatively-oriented, local orthonormal frame of $E$. One can check that $\tau_q=-\segno$.  
Let us make the following change of coordinates
$$\tilde{x}=x, ~~~~~\tilde{y}=\segno(y-x^2\psi(x)).$$
In these new coordinates, $X$ and $Y$ become
$$X=(1,-\segno(2\tilde{x}\psi(\tilde{x})+\tilde{x}^2\psi'(\tilde{x}))), ~~~~~Y=(0,\tilde{y} e^{\xi(\tilde{x},\segno\tilde{y}+\tilde{x}^2\psi(\tilde{x}))})$$
and $W\cap U$ is the $\tilde{x}$-axis. In the following, to simplify notations, we  omit the tildes  and denote the function $\xi(\tilde{x},\alpha\tilde{y}+\tilde{x}^2\psi(\tilde{x})))$ by $\xi(x,y)$.
Since $\frac{\partial\psi}{\partial \th}(\th,0)$ is tangent to $W$, by the Preparation Theorem  \cite{malgrange}
 there exist 
$h_1:\R\rightarrow\R$, $h_2,h_3:\R^2\rightarrow\R$ smooth functions such that $h_1(x)\neq 0$ for every $x\in\R$ and for $\Psi(\th,s)\in  U$
$$
\frac{\partial\psi}{\partial \th}(\th,s)=(h_1(x)+y h_2(x,y), y h_3(x,y))
$$
where $(x,y)$ are the coordinates of the point $\Psi(\th,s)$.
This implies that
\begin{eqnarray*}
\frac{\partial\psi}{\partial
\th}(\th,s)&=&(h_1(x)+yh_2(x,y))X(x,y)+\\
&+&\frac{yh_3(x,y)+\segno(h_1(x)+yh_2(x,y))(2x\psi(x)+x^2\psi'(x))}{y e^{\xi(x,y)}}Y(x,y).
\end{eqnarray*}
Let $(\zeta,\rho)$ be the local orthonormal frame of $E$ such that $X=\f\circ\zeta$ and $Y=\f\circ\rho$. 
From equation~\r{sez}, it follows that
\begin{eqnarray*}
\sigma(x,y)&=&-\segno\frac{\sign(y)}{l(x,y)}\frac{yh_3(x,y)+\segno(h_1(x)+yh_2(x,y))(2x\psi(x)+x^2\psi'(x))}{y e^{\xi(x,y)}}\zeta+\\
&+&\segno\frac{\sign(y)}{l(x,y)}(h_1(x)+yh_2(x,y))\rho,
\end{eqnarray*}
where
$$
l(x,y)=\sqrt{(h_1(x)+yh_2(x,y))^2+\left(\frac{yh_3(x,y)+\segno(h_1(x)+yh_2(x,y))(2x\psi(x)+x^2\psi'(x))}{y e^{\xi(x,y)}}\right)^2}.
$$
Notice that for $x=0, y\neq 0$ we have 
$$
\sigma(0,y)=\frac{\segno\ \sign(y)}{\sqrt{(h_1(0)+yh_2(0,y))^2+h_3(0,y)^2e^{-2\xi(0,y)}}}(-e^{-\xi(0,y)}h_3(0,y)\zeta+(h_1(0)+yh_2(0,y))\rho),
$$
whence the limit of $\sigma$ as $(x,y)$ tends to $(0,0)$ does not exist.
Let us compute the index of $\sigma$ at $q=(0,0)$. 
Using Taylor expansions of the components of $\sigma$ in the basis $(\zeta,\rho)$ we find
\begin{equation}\label{expan}
\sigma(x,y)=\segno\frac{\sign(y)}{l(x,y)}\left(-\left(\frac{yh_3(0,0)+2\segno xh_1(0)\psi(0)}{y e^{\xi(0,0)}}+O(\sqrt{x^2+y^2})\right)\zeta+(h_1(0)+O(\sqrt{x^2+y^2}))\rho\right).
\end{equation}
Take a circle $t\mapsto(r \cos t,r\sin t)$ of radius $r$ centered at $(0,0)$ and assume $r$  so small that $(0,0)$ is the unique singularity of $\sigma$ on the closed disk of radius $r$. By definition, $i((0,0),\sigma)$ is half the  degree of the map from the circle $S^r$ to $\R/\pi\Z$ that associates to each point the angle  between $\mathrm{span}(\sigma)$ and the $\zeta$. Using~\r{expan},  this angle is 
$$
a(x,y)=-\arctan{\left(\frac{yh_1(0)e^{\xi(0,0)}}{yh_3(0,0)+2\segno x\psi(0)h_1(0)}+O(\sqrt{x^2+y^2})\right)}.
$$
Computing $a$ along the curve $x(t)=r\cos t,\,y(t)=r\sin t$, we find
$$
a(r\cos t,r\sin t)=-\arctan{\left(\frac{\sin(t) h_1(0)e^{\xi(0,0)}}{\sin(t) h_3(0,0)+2\segno\cos(t)\psi(0)h_1(0)}+O(r)\right)}.
$$
Hence, by letting $r$ go to zero, we are left to compute the degree of the map $\tilde{a}:[0,2\pi)\rightarrow[0,\pi)$ where
$$
\tilde{a}(t)=-\arctan{\left(\frac{\sin(t) h_1(0)e^{\xi(0,0)}}{\sin(t) h_3(0,0)+2\segno\cos(t)\psi(0)h_1(0)}\right)}.
$$
 Since zero is a regular value of $\tilde{a}$, the degree of $\tilde{a}$ is 
$$
\sum_{t\in[0,2\pi)\mid\tilde{a}(t)=0}\sign(\tilde{a}'(t))=\sign(\tilde{a}'(0))+\sign(\tilde{a}'(\pi))=-2\segno,
$$
where the last equality follows from $\tilde{a}'(0)=\tilde{a}'(\pi)=-\segno e^{\xi(0,0)}/(2\psi(0))$. Hence, $i(q,\sigma)=-\segno$. Since $\tau_q=-\segno$, the lemma is proved  (see Remark~\ref{nuovopn}).
\hfill$\blacksquare$

 Let $C(\Zz)$ denote the set of connected component of $\Zz$. Let  $\tilde{\Zz}=\coprod_{W\in C(\Zz)}S^1$  and consider an orientation-preserving diffeomorphism $\Psi:\tilde{\Zz}\times(-1,1)\rightarrow \coprod_{W\in C(\Zz)}\cil$ such that $\Psi|_{\tilde{\Zz}\times \{0\}}$ is an orientation-preserving diffeomorphism onto $\Zz$.
Applying Lemma~\ref{unico} to every  $W\in C(\Zz)$ and reducing, if necessary, the cylinders $\cil$, we can assume  that 
the set of singularities of $\sigma$ on $\CIL=\coprod_{W\in C(\Zz)}\cil$ is ${\cal T}$. Then $\sigma:\CIL\setminus{\cal T}\rightarrow SE$ is  continuous. Moreover, by equation~\r{verita},
$$
\sum_{q\in{\cal T}}i(q,\sigma)=\tau({\cal S}).
$$
 Extend $\sigma$ to $M\setminus \CIL$. By a transversality argument, we can assume that the extended section has only isolated singularities $\{p_1,\dots, p_k\}\in M\setminus\Zz$. Since 
$$
\E(E)=\sum_{j=1}^ki(p_j,\sigma)+\sum_{q\in{\cal T}}i(q,\sigma)= \sum_{j=1}^ki(p_j,\sigma)+\tau({\cal S}),
$$
we are left to prove that 
\begin{equation}\label{claim}
 \sum_{j=1}^ki(p_j,\sigma)=\chi(M^+)-\chi(M^-).
\end{equation}

To this aim, consider the vector field $F=\f \circ\sigma$. $F$ satisfies $G(F,F)\equiv 1$, where $G(\cdot,\cdot)$ is defined as in Remark~\ref{equiv} and the set of singularities of $F|_{M\setminus\Zz}$ is exactly $\{p_1,\dots,p_k\}$. Let us compute the index of $F$ at a singularity $p\in \{p_1,\dots, p_k\}$. Since $\f:E|_{M^+}\rightarrow TM^+$ preserves the orientation and $\f:E|_{M^-}\rightarrow TM^-$ reverses the orientation, it follows that $i(p,F)=\pm i(p,\sigma)$, if $p\in M^\pm$. 
Therefore,
\begin{equation}\label{five}
 \sum_{j=1}^ki(p_j,\sigma)=\sum_{j\mid p_j\in M^+}i(p_j,F)-\sum_{j\mid p_j\in M^-}i(p_j,F).
\end{equation}
The theorem is proved if we show that
\begin{equation}\label{fine}
\sum_{j\mid p_j\in M^+}i(p_j,F)=\chi(M^+),\,\,\sum_{j\mid p_j\in M^-}i(p_j,F)=\chi(M^-).
\end{equation}
To deduce equation~\r{fine}, define $N^+=M^+\setminus \Psi(\tilde{\Zz}\times(0,1/2))$. Notice that, by construction, $\sigma|_{\Psi(\tilde{\Zz}\times\{1/2\})}$ is non-singular, hence the same is true for $F|_{\Psi(\tilde{\Zz}\times\{1/2\})}$.
Moreover, the almost-Riemannian angle  between $T_q(\partial N^+)$ and  $\mathrm{span }(F(q))$ is constantly equal to $\pi/2$. Hence $F|_{\partial N^+}$ points towards $N^+$ and applying the Hopf's Index Formula to every connected component of $N^+$ we conclude that 
$$
\sum_{j\mid p_j\in M^+}i(p_j,F)=\sum_{j\mid p_j\in N^+}i(p_j,F)=\chi(N^+)=\chi(M^+).
$$ 
Similarly, we find 
$$
 \sum_{j\mid p_j\in M^-}i(p_j,F)=\chi(M^-).
$$\hfill$\blacksquare$

\section{$\cal S$-integrability in presence of tangency points}\label{integr}

\subsection{Numerical simulations}\label{numsim}

In this section we provide some numerical simulations hinting that, when ${\cal T}\neq\emptyset$,
$$
\int_{M_\eps} K dA_s
$$
does not converge, in general, as $\eps$ tends to zero. 

From the proof of Theorem~\ref{gbsm} we know that far from tangency points the integral of the geodesic curvature along $\partial M^+_\eps$ and $\partial M^-_\eps$ offset each other for $\eps$ going to zero. Hence, to understand whether the presence of a tangency point may lead to  non-$\cal S$-integrability of $K$ it is sufficient to compute the geodesic curvature of $\partial M^+_\eps$ and   $\partial M^-_\eps$ in a neighborhood of such a point. More precisely consider the almost-Riemannian structure $(E,\f,\langle\cdot,\cdot\rangle)$ on $M=\R^2$ for which $E=\R^2\times \R^2$, $\f((x,y),(a,b))=((x,y),(a,b(y-x^2)))$ and $\langle\cdot,\cdot\rangle$ is the canonical scalar product. For this system one has
\bqn
K=\frac{-2\,\left( 3\,x^2 + y \right)}
  {{\left( x^2 - y \right) }^2}.
\eqnn
The graph of $K$ is illustrated in
Figure~\ref{curvature-tang}.
Notice that
$\limsup_{q\to (0,0)}
K(q)=+\infty$ and\linebreak $\liminf_{q\to (0,0)} K(q)=-\infty$. This situation is different from the Grushin case where $K(q)$ diverges to $-\infty$ as $q$ approaches $\Zz$. 

For every $\eps>0$, the sets $\partial M^+_\eps$ and   $\partial M^-_\eps$ are smooth manifolds except at their intersections with the vertical axis $x=0$, which is the cut locus for the problem of minimizing the distance from $\Zz=\{(x,x^2)\mid x\in\R \}$. 
Fix $0<a<1$ and consider the two geodesics starting from the point $(a,a^2)$ and minimizing (locally) the distance from  $\Zz$. Let $P^+$ and $P^-$ be the two points along these geodesics at distance $\eps$ from $\Zz$. Denote by   $\gamma^+$ and $\gamma^-$ the portions of  $\partial M^+_\eps$ and   $\partial M^-_\eps$ connecting the vertical axis to the points $P^+$ and $P^-$, oriented as in Figure~\ref{f-numeriKa}. It is easy to approximate numerically $\g^+$ and $\g^-$ by broken lines, but
the evaluation of the integral of their geodesic curvatures is very unstable since its computation involves the second derivative of the curve parameterized by arclength. To avoid this problem, we rather apply the Riemannian Gauss--Bonnet formula 
on the regions $\Omega^+$ and $\Omega^-$ introduced in Figure~\ref{f-numeriKa}. This works better since the integral of the Gaussian curvature on $\Omega^+$ and $\Omega^-$ is numerically stable, and the integral of the geodesic curvature on horizontal and vertical segments can be  computed analytically (in particular it is always zero on horizontal segments). Figure~\ref{gggpunti-bis}  shows the value of 
$$\eps \left(\int_{\gamma^+} K_g ds - \int_{\gamma^-} K_g ds   \right)$$ for $a=0.1$ and $\eps$ varying in the interval $[0.01,0.04]$. The graph seems to converge as $\eps$ tends to zero to a nonzero constant,  strongly hinting at the divergence of  
$\int_{M_\eps} K dA_s$. 

\ppotR{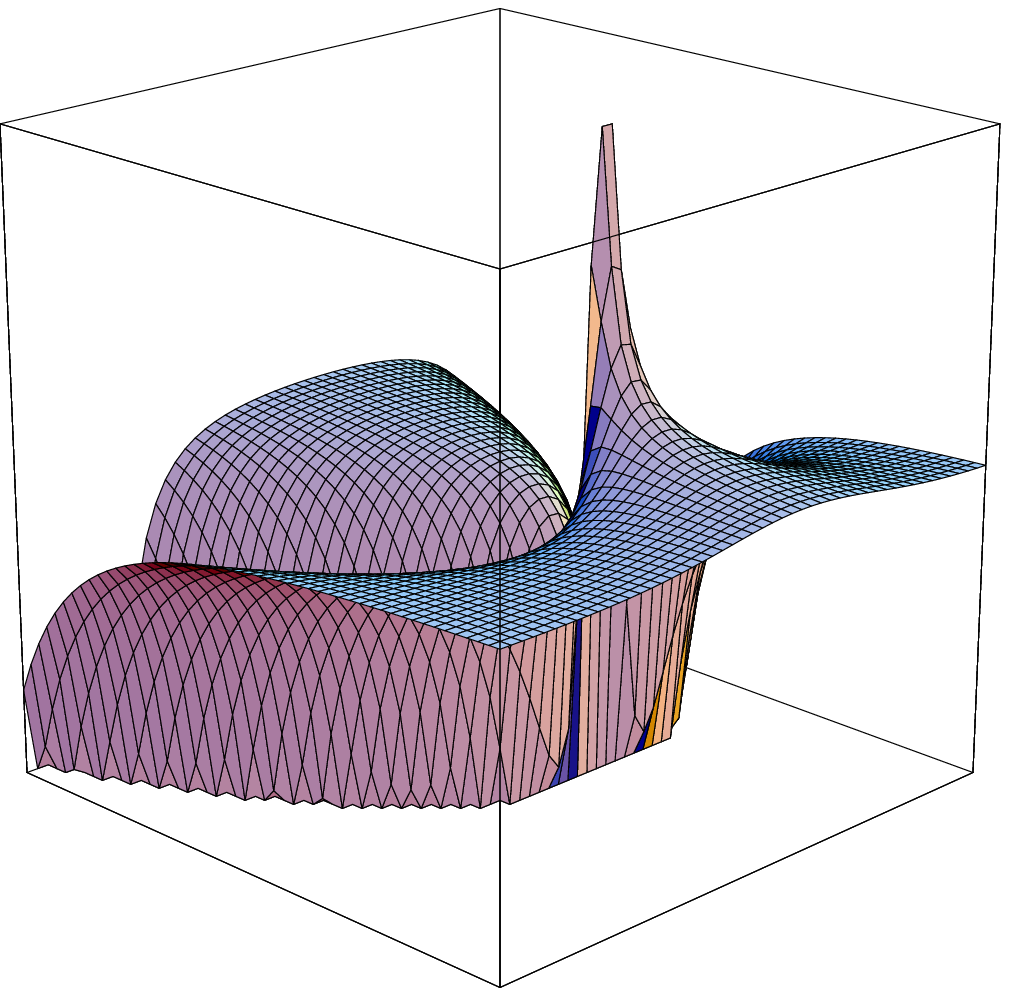}{Graph of $K$ for $\Delta = \mathrm{span}((1,0),(0,y-x^2))$}{7}
\ppotR{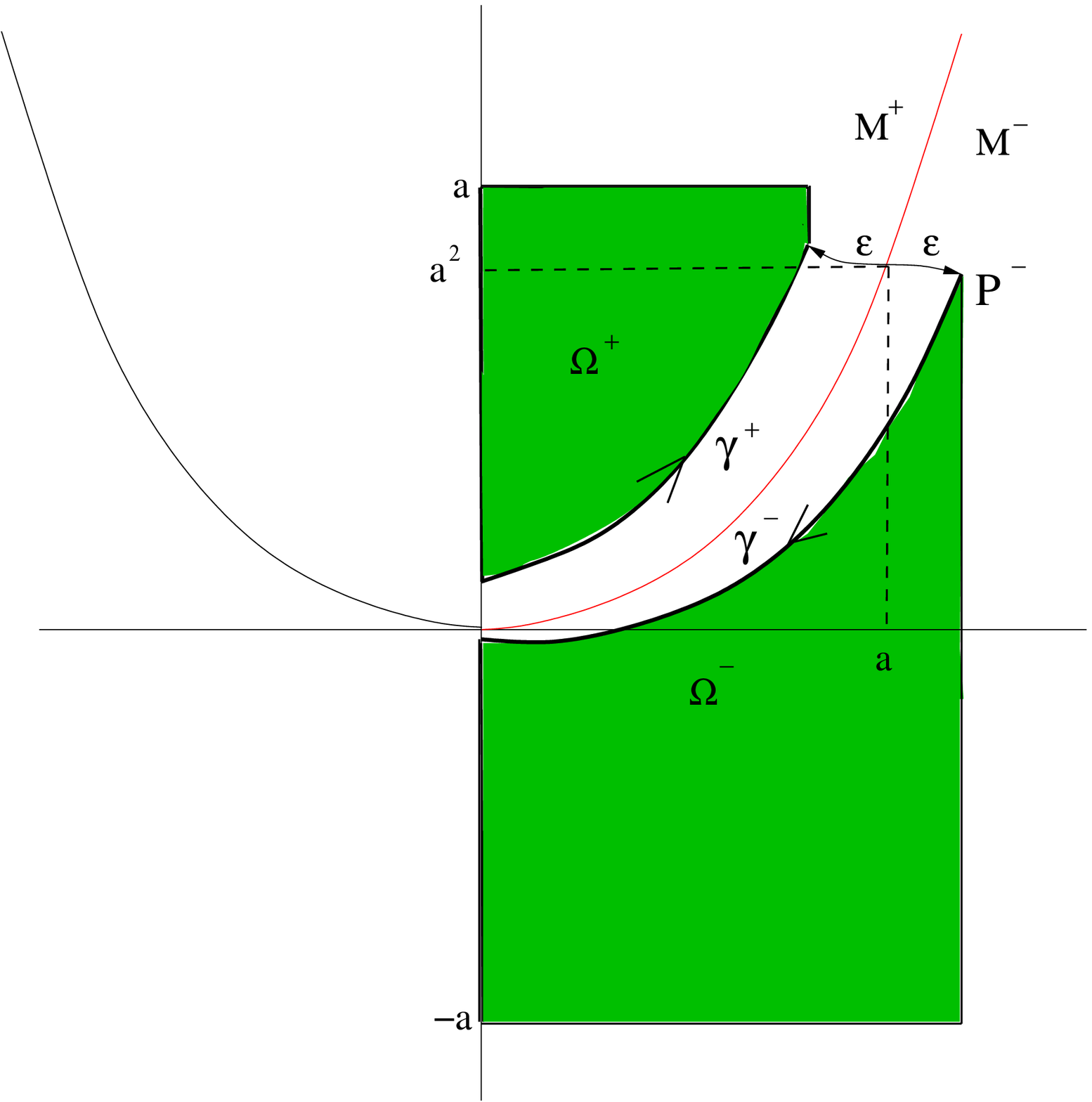}{Regions $\Omega^\pm$ where to apply Riemannian Gauss--Bonnet formula}{9}
\ppotR{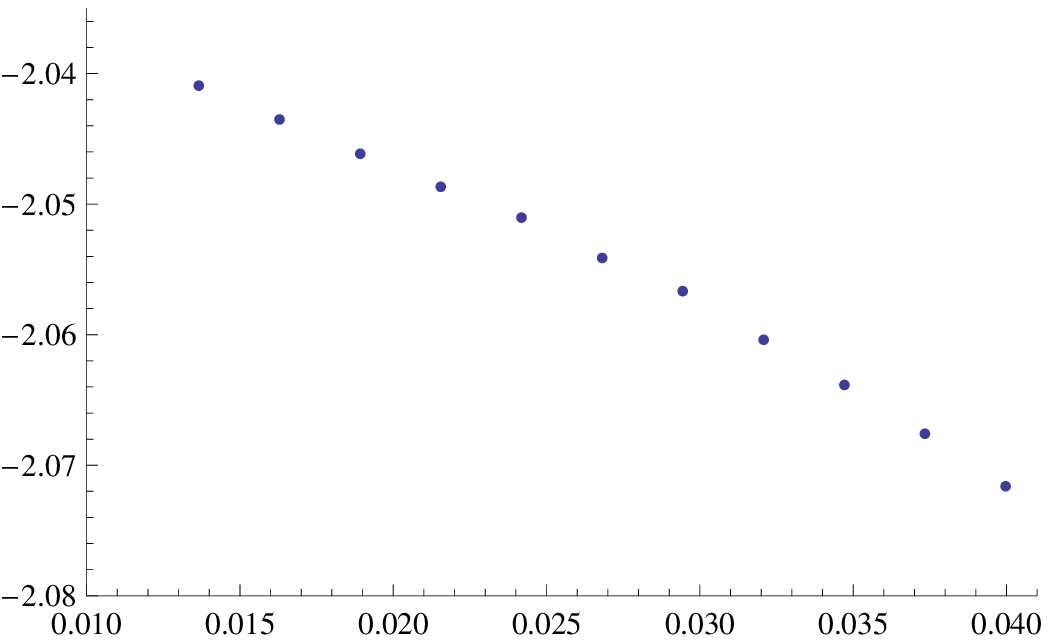}{Divergence of the $\cal S$-integral of $K$}{10}
\subsection{More general notion of $\cal S$-integrability}\label{3int}

The simulations of the previous section lead us to 
introduce the following alternative notion of integrability. 

\bdeff {\bf (3-scale ${\cal S}$-integrability)}
  Let $q\in{\cal T}$ and  $U^q$ be a neighborhood of 
  $q$ such that an orthonormal frame for $G$ on $U^q$ is given by the normal form (F3). For $\delta_1$, $\delta_2>0$ sufficiently small  the rectangle $[-\delta_1,\delta_1]\times[-\delta_2,\delta_2]$ is a subset of $U^q$ denoted by $B^q_{\delta_1,\delta_2}$.
For every $\eps>0$, define
\bqn
M_{\eps,\delta_1,\delta_2}=M_\eps \setminus  \bigcup_{q\in{\cal T} }  B^q_{\delta_1,\delta_2}. 
\eqn
We say that $K$ is 3-scale ${\cal S}$-integrable 
if 
\bqn\label{box}
\lim_{\delta_1\to0}\lim_{\delta_2\to0}\lim_{\eps\to0}\int_{M_{\eps,\delta_1,\delta_2}} K dA_s
\eqn
exists, is finite and does not depend on the choice of the normal form. In this  case we denote such limit by 
$\treint_M K dA_s$.
\edeff

\brem
Notice that if ${\cal T}=\emptyset$, then the concepts of ${\mathcal S}$-integrability and 3-scale ${\mathcal S}$-integrability coincide. 

The order in which the limits are taken in~\r{box} is important. Indeed, if the order is permuted, then the result given in Theorem~\ref{t-GB} does not hold anymore.

Recall that the normal form (F3) is not totally intrinsic, since the functions $\psi$ and $\phi$ depend on the choice of a parametrized smooth curve passing through the tangency point  and transversal to the distribution at the point. 
\erem


\subsection{Proof of Theorem~\ref{t-GB}}\label{secpro}
Let us recall the following  Gauss--Bonnet-like formula for domains whose boundary is $\con^2$ in  a neighborhood  of $\Zz$.
\bt[\cite{high-order}, Theorem 5.2]
\label{t-boundary}
Let $U$ be 
an open bounded connected subset 
of $M$ such that
{\bf i)} $\overline U$ contains only ordinary and Grushin points, 
{\bf ii)}  $\partial U$ is piecewise $\con^2$,
{\bf iii)}  $\partial U$ is  $\con^2$ in a neighborhood of $\Zz$,
{\bf iv)}  $\partial U$ is the union of the supports
of a finite set of curves $\g^1,\dots,\g^m$  that are admissible for $\Delta$ and of finite length.  

Define $U^\pm_\eps=M^\pm_\eps\cap U$. Then the following limits exist and are finite
\bqn
\int_UK dA_s&:=&\lim_{\eps\to 0}\int_{U^+_\eps\cup U^-_\eps}K dA_s,\label{Ku}\\
\int_{\partial U}k_gd\sigma_s&:=&\lim_{\eps\to 0}\left(\int_{\partial U\cap \partial U^+_\eps}k_g d\sigma- 
\int_{\partial U\cap \partial U^-_\eps}k_g d\sigma\right),\label{kgu}
\eqn
where we interpret each 
integral $\int_{\partial U\cap \partial U^\pm_\eps}k_g d\sigma$ as the sum of the  integrals  along the 
$\con^2$ portions of $\partial U\cap \partial U^\pm_\eps$, plus the sum of the angles at the points of $\partial U\cap \partial U^\pm_\eps$ where 
$\partial U$ is not $\con^1$. Moreover, we have
\bqn
\int_UK dA_s+\int_{\partial U}k_gd\sigma_s=2\pi(\chi(U^+)-\chi(U^-)).
\eqn
\et

Fix $\delta_1$ and $\delta_2$ in such a way that  the rectangles 
$B^q_{\delta_1,\delta_2}$ are pairwise disjoint and
$\Zz\cap\partial B^q_{\delta_1,\delta_2}\subset [-\delta_1,\delta_1]\times\{\delta_2\}$, for every $q\in{\cal T}$.
By construction, 
$\partial B^q_{\delta_1,\delta_2}$ is admissible and has finite length for every $q\in {\cal T}$. Hence we can take  $M_\eps\setminus \bigcup_{q\in {\cal T}} B^q_{\delta_1,\delta_2}$
 as $U$ in Theorem~\ref{t-boundary}.
As a consequence, we have
 \bqn
\lim_{\eps\to0}\int_{M_{\eps,\delta_1,\delta_2}} K dA_s+\sum_{q\in{\cal T}}\int_{\partial B^q_{\delta_1,\delta_2}}k_g d\sigma_s&=&2\pi(\chi(M^+\setminus \bigcup_{q\in{\cal T}} B^q_{\delta_1,\delta_2})-  \chi( M^-\setminus \bigcup_{q\in{\cal T}} B^q_{\delta_1,\delta_2} ))\nn\\
&=&2\pi(\chi(M^+)-\chi(M^-))=2\pi(\E(E)-\tau({\cal S})),\nn
\eqn
where the last equality follows from Theorem~\ref{converse}.
We are left to prove that, for a fixed $q\in{\cal T}$, 
$$
\lim_{\delta_1\to 0}\lim_{\delta_2\to 0}\int_{\partial B^q_{\delta_1,\delta_2}}k_g d\sigma_s=-2\pi \tau_{q} 
$$
(see Remark~\ref{nuovopn}).
 In order to prove it, let us work with the normal form $(F3)$ and assume that $M^+ \cap U^q$ is the set $\{y-x^2\psi(x)< 0\}\cap U^q$, the proof for the opposite situation being analogous. On one hand, one can check that $\tau_{q}=1$. On the other hand, the geodesic curvature along $[-\delta_1,\delta_1]\times\{\delta_2\}$ and along $[-\delta_1,\delta_1]\times\{-\delta_2\}$ is zero, the two segments being 
the support of geodesics. Hence
$$
\int_{\partial B^q_{\delta_1,\delta_2}}k_g d\sigma_s=\int_{\{\delta_1\}\times[-\delta_2,\delta_2]}k_g d\sigma_s+
\int_{\{-\delta_1\}\times[-\delta_2,\delta_2]}k_g d\sigma_s+ \sum_{j=1}^4\alpha_j
$$ where the last term is the sum of the values of the angles of the box and is equal to $-2\pi$. Indeed, because of the diagonal form of the metric with respect to the chosen coordinates, each angle has value $-\frac \pi 2$. The first two terms are well defined and tend  to zero when $\delta_2$ tends to zero. Hence 
$$
\lim_{\delta_1\to 0}\lim_{\delta_2\to 0}\int_{\partial B^q_{\delta_1,\delta_2}}k_g d\sigma_s=-2\pi=-2\pi \tau_{q}.
$$\hfill$\blacksquare$

\bibliographystyle{abbrv}
\bibliography{biblio_euler}

\end{document}

%% file: boxtrasv.pstex_t
\begin{picture}(0,0)%
\includegraphics{boxtrasv.pstex}%
\end{picture}%
\setlength{\unitlength}{1381sp}%
\begingroup\makeatletter\ifx\SetFigFont\undefined%
\gdef\SetFigFont#1#2#3#4#5{%
  \reset@font\fontsize{#1}{#2pt}%
  \fontfamily{#3}\fontseries{#4}\fontshape{#5}%
  \selectfont}%
\fi\endgroup%
\begin{picture}(13524,9849)(-311,-9448)
\put(3076,-511){\makebox(0,0)[lb]{\smash{{\SetFigFont{9}{10.8}{\rmdefault}{\mddefault}{\updefault}{\color[rgb]{0,0,0}$\bD$}%
}}}}
\put(6151,-8611){\makebox(0,0)[lb]{\smash{{\SetFigFont{9}{10.8}{\rmdefault}{\mddefault}{\updefault}{\color[rgb]{0,0,0}$\gamma_{-\delta_2}(\delta_1)$}%
}}}}
\put(5401,-5611){\makebox(0,0)[lb]{\smash{{\SetFigFont{9}{10.8}{\rmdefault}{\mddefault}{\updefault}{\color[rgb]{0,0,0}$q$}%
}}}}
\put(9751,-2011){\makebox(0,0)[lb]{\smash{{\SetFigFont{9}{10.8}{\rmdefault}{\mddefault}{\updefault}{\color[rgb]{0,0,0}$M^+$}%
}}}}
\put(11551,-4036){\makebox(0,0)[lb]{\smash{{\SetFigFont{9}{10.8}{\rmdefault}{\mddefault}{\updefault}{\color[rgb]{0,0,0}$M^-$}%
}}}}
\put(1276,-211){\makebox(0,0)[lb]{\smash{{\SetFigFont{9}{10.8}{\rmdefault}{\mddefault}{\updefault}{\color[rgb]{0,0,0}$\Zz$}%
}}}}
\put(10276,-6811){\makebox(0,0)[lb]{\smash{{\SetFigFont{9}{10.8}{\rmdefault}{\mddefault}{\updefault}{\color[rgb]{0,0,0}$\gamma_{\delta_2}(\delta_1)$}%
}}}}
\put(5776,-3586){\makebox(0,0)[lb]{\smash{{\SetFigFont{9}{10.8}{\rmdefault}{\mddefault}{\updefault}{\color[rgb]{0,0,0}$w(\delta_2)$}%
}}}}
\put(4276,-7186){\makebox(0,0)[lb]{\smash{{\SetFigFont{9}{10.8}{\rmdefault}{\mddefault}{\updefault}{\color[rgb]{0,0,0}$w(-\delta_2)$}%
}}}}
\put(376,-6436){\makebox(0,0)[lb]{\smash{{\SetFigFont{9}{10.8}{\rmdefault}{\mddefault}{\updefault}{\color[rgb]{0,0,0}$\gamma_{-\delta_2}(-\delta_1)$}%
}}}}
\put( 76,-2836){\makebox(0,0)[lb]{\smash{{\SetFigFont{9}{10.8}{\rmdefault}{\mddefault}{\updefault}{\color[rgb]{0,0,0}$\gamma_{\delta_2}(-\delta_1)$}%
}}}}
\end{picture}%

%% file: tauq.pstex_t
\begin{picture}(0,0)%
\includegraphics{tauq.pstex}%
\end{picture}%
\setlength{\unitlength}{1579sp}%
\begingroup\makeatletter\ifx\SetFigFont\undefined%
\gdef\SetFigFont#1#2#3#4#5{%
  \reset@font\fontsize{#1}{#2pt}%
  \fontfamily{#3}\fontseries{#4}\fontshape{#5}%
  \selectfont}%
\fi\endgroup%
\begin{picture}(14049,4374)(-236,-6523)
\put( 76,-4336){\makebox(0,0)[lb]{\smash{{\SetFigFont{8}{9.6}{\rmdefault}{\mddefault}{\updefault}{\color[rgb]{0,0,0}$W$}%
}}}}
\put(5026,-4036){\makebox(0,0)[lb]{\smash{{\SetFigFont{8}{9.6}{\rmdefault}{\mddefault}{\updefault}{\color[rgb]{0,0,0}$M^+$}%
}}}}
\put(12376,-4711){\makebox(0,0)[lb]{\smash{{\SetFigFont{8}{9.6}{\rmdefault}{\mddefault}{\updefault}{\color[rgb]{0,0,0}$M^-$}%
}}}}
\put(4951,-4786){\makebox(0,0)[lb]{\smash{{\SetFigFont{8}{9.6}{\rmdefault}{\mddefault}{\updefault}{\color[rgb]{0,0,0}$M^-$}%
}}}}
\put(7351,-4336){\makebox(0,0)[lb]{\smash{{\SetFigFont{8}{9.6}{\rmdefault}{\mddefault}{\updefault}{\color[rgb]{0,0,0}$W$}%
}}}}
\put(12526,-4036){\makebox(0,0)[lb]{\smash{{\SetFigFont{8}{9.6}{\rmdefault}{\mddefault}{\updefault}{\color[rgb]{0,0,0}$M^+$}%
}}}}
\put(2551,-5611){\makebox(0,0)[lb]{\smash{{\SetFigFont{8}{9.6}{\rmdefault}{\mddefault}{\updefault}{\color[rgb]{0,0,0}$\tau_q=-1$}%
}}}}
\put(10051,-5611){\makebox(0,0)[lb]{\smash{{\SetFigFont{8}{9.6}{\rmdefault}{\mddefault}{\updefault}{\color[rgb]{0,0,0}$\tau_q=1$}%
}}}}
\put(10426,-4861){\makebox(0,0)[lb]{\smash{{\SetFigFont{6}{7.2}{\rmdefault}{\mddefault}{\updefault}{\color[rgb]{0,0,0}$q$}%
}}}}
\put(2851,-4711){\makebox(0,0)[lb]{\smash{{\SetFigFont{6}{7.2}{\rmdefault}{\mddefault}{\updefault}{\color[rgb]{0,0,0}$q$}%
}}}}
\put(1351,-2761){\makebox(0,0)[lb]{\smash{{\SetFigFont{8}{9.6}{\rmdefault}{\mddefault}{\updefault}{\color[rgb]{0,0,0}$\bD$}%
}}}}
\put(8776,-2761){\makebox(0,0)[lb]{\smash{{\SetFigFont{8}{9.6}{\rmdefault}{\mddefault}{\updefault}{\color[rgb]{0,0,0}$\bD$}%
}}}}
\end{picture}%